\documentclass[12pt,reqno]{amsart}
\usepackage{}
\usepackage{mathrsfs}
\usepackage{amssymb}
\usepackage{amsthm}
\usepackage{mathrsfs}
\usepackage[centertags]{amsmath}
\usepackage{amsfonts}
\usepackage[numbers,sort&compress]{natbib}
\usepackage{color}
\usepackage{extarrows}
\usepackage{fullpage}
\usepackage[colorlinks=true,  linkcolor=blue, citecolor=blue]{hyperref}
\input amssym.def
\input amssym.tex

\date{}

\newtheorem{Theorem}{Theorem}[section]

\newtheorem{Lemma}{Lemma}[section]

\newcommand\R{\mbox{\bf R}}

\newcommand\Z{\mbox{\bf Z}}
\newcommand\z{\mbox{\bf z}}
\newcommand\SR{\mbox{\scriptsize\bf R}}

\newcommand{\definition}{{\lower .5ex
  \hbox{$\>\>\stackrel{\triangle}{=}\>\>$} }}
\newcommand\supp{\mathop{\rm supp}}



\begin{document}

\baselineskip=22pt
\thispagestyle{empty}

\begin{center}
{\Large \bf  The probabilistic convergence
 problem of density functions related to $\partial_{x}^{3}+\partial_{x}^{-1}$}\\[1ex]

{ Xiangqian Yan\footnote{Email: yanxiangqian213@126.com}$^{a}$,\,
Yongsheng Li\footnote{Email:  yshli@scut.edu.cn}$^{a}$,\,Wei Yan\footnote{Email:  011133@htu.edu.cn}$^{b*}$}\\[1ex]

{$^a$School of Mathematics,
 South China University of Technology,}\\
{ Guangzhou, Guangdong 510640, P. R. China}\\[1ex]

{$^b$School of Mathematics and Statistics, Henan
Normal University,}\\
{Xinxiang, Henan 453007, P. R. China}\\[1ex]

\end{center}
\noindent{\bf Abstract.}
In this article, by  using full randomization introduced by   Hadama and  Yamamoto (Probabilistic Strichartz estimates in Schatten classes
and their applications to Hartree equation,  J. Math. Phys. 67(2026), 35pp) and high-low frequency technique as well as the property of $\mathfrak{S}^{2}$,
      we establish the probabilistic convergence of the density function
       related to $\partial_{x}^{3}+\partial_{x}^{-1}$
       on $\R$, which extends the Theorem 1.3 of Yan et al.
          (Convergence problem of
Ostrovsky equation with rough data and random data,  Indiana Univ. Math. J. 71(2022), 1897-1921.).

 \bigskip

\medskip

\medskip
\noindent {\bf Corresponding Author:} Wei Yan

\medskip
\noindent {\bf Email Address:}  011133@htu.edu.cn

\medskip
\noindent {\bf MSC2020-Mathematics Subject Classification}: Primary-35Q40, 47B10; Secondary-42B37

\leftskip 0 true cm \rightskip 0 true cm

\newpage

\baselineskip=20pt

\bigskip
\bigskip
\tableofcontents

\section{Introduction}
\bigskip

\setcounter{Theorem}{0} \setcounter{Lemma}{0}\setcounter{Definition}{0}\setcounter{Proposition}{1}

\setcounter{section}{1}
\subsection{Research background and current status}
In this paper, the following infinite systems of linear Ostrovsky equations are studied
\begin{eqnarray}
&&\frac{du_{j}}{dt}=-(\partial_{x}^{3}+\partial_{x}^{-1})u_{j},
\,j\in \mathbb{N}^{+}, x\in\R,\,t\in I\subset\R,\label{1.01}\\
&&u_{j}(0,x)=f_{j}(x),\label{1.02}
\end{eqnarray}
where $(f_j)_{j=1}^{\infty}$ is an orthonormal system in $L^{2}(\R)$.
From \eqref{1.01},
 we can get its the operator valued form
\begin{eqnarray}
&&\frac{d\gamma(t)}{dt}=\left[-(\partial_{x}^{3}+\partial_{x}^{-1}), \gamma\right], \label{1.03}\\
&&\gamma(0)=\gamma_{0}=\sum_{j=1}^{\infty}\lambda_{j}|f_{j}\rangle\langle f_{j}|.\label{1.04}
\end{eqnarray}
Evidently,
\begin{eqnarray}
&\gamma(t):=e^{-t(\partial_{x}^{3}+\partial_{x}^{-1})}\gamma_{0}
e^{t(\partial_{x}^{3}+\partial_{x}^{-1})}=\sum\limits_{j=1}^{\infty}\lambda_{j}
\left|e^{-t(\partial_{x}^{3}+\partial_{x}^{-1})}f_{j}\right\rangle\left\langle
 e^{-t(\partial_{x}^{3}+\partial_{x}^{-1})}f_{j}\right|\label{1.05}
\end{eqnarray}
is the solution to (\ref{1.03})-(\ref{1.04}).
Here $\lambda_{j}$ is the eigenvalue of the operator $\gamma_{0}$, and
 $|u\rangle\langle v| f\rightarrow\langle v,f\rangle u=u\int_{I}\overline{v}f dx$,
 $I\subset\R$.
From \eqref{1.05}, we derive the kernel of $\gamma(t)$ as follows:
\begin{eqnarray}
&K(x,y,t)=\sum\limits_{j=1}^{\infty}\lambda_{j}
\left(e^{-t(\partial_{x}^{3}+\partial_{x}^{-1})}f_{j}\right)(t,x)\overline
{\left(e^{-t(\partial_{x}^{3}+\partial_{x}^{-1})}f_{j}\right)}(t,y).\label{1.06}
\end{eqnarray}
By using \eqref{1.06}, we have the density function of  $\gamma(t)$ as follows:
\begin{eqnarray*}
\rho_{\gamma(t)}
=K(x,x,t)=\sum\limits_{j=1}^{\infty}\lambda_{j}
\left|e^{-t(\partial_{x}^{3}+\partial_{x}^{-1})}f_{j}\right|^{2}.
\end{eqnarray*}

\subsection{The background of pointwise convergence and probabilistic convergence for single equation}

The  pointwise convergence problem of Schr\"odinger flow  is originally introduced by Carleson \cite{C1979}.
In the past few decades, some people have studied the pointwise convergence problem
in high-dimensional situations, which can be seen in
\cite{B1995,B2016,CLV2012,C1982,DK1982,DG,D2017,DGL2017,DZ2019,DGLZ2018,L2006,
 LR,LR2017,MV2008,S,V1988,Z2014}. Except the endpoint, Du and her co-authors \cite{DGL2017,DZ2019}
 obtained optimal results for the  pointwise convergence problem of Schr\"odinger flow
 in all dimensions.
 Recently, Compaan et al. \cite{CLS2021} applied the  technique of initial data randomization proposed by
Lebowitz-Rose-Speer  \cite{LRS1988} and  Bourgain  \cite{B1994,B1996}  and later developed by Burq-Tzvetkov
  \cite{BT2007,BT2008I,BT2008II} to study the  pointwise convergence
   of the Schr\"{o}dinger flow with random data.
    Moreover, by using the  technique of initial data randomization, Wang et al. \cite{WYY2021} established the pointwise
 convergence for the Schr\"odinger equation on $\mathbf{T}^{d}$, and the unit ball $\Theta=\{x\in\R^{3}:|x|\leq 1\}$ with random data.
  By using  maximal-in-time estimates, the density theorem in Fourier-Lebesgue spaces and the  technique of initial data randomization,
    Yan et al. \cite{YZY2022} established  the pointwise
 convergence of the Schr\"odinger equation in Fourier-Lebesgue spaces with rough data and random data.
      Very recently, by using high-low frequency method and the  technique of initial data randomization, Yan et al. \cite{YZDY2022} proved  the pointwise convergence
   of the Ostrovsky equation with rough data and random data.

\subsection{The background of pointwise convergence and probabilistic convergence for the density function of
operators}

When studying pointwise convergence for linear equations, a standard method is to consider the maximal-in-time estimates associated with the solution.
For the one-dimensional linear Schr\"{o}dinger equation
\begin{eqnarray*}
&&iu_{t}+\partial^{2}_{x}u=0,\, (x,t)\in\R\times\R,\\
&&u(0,x)=f(x),
\end{eqnarray*}
Kenig et al. \cite{KPV1991} established the following maximal-in-time estimate
\begin{eqnarray}
&&\left\|\sup_{t\in\SR}e^{it\partial_{x}^{2}}f\right\|_{L_{x}^{4}}\leq C\|f\|_{\dot{H}^{\frac{1}{4}}}.\label{1.07}
\end{eqnarray}
By using \eqref{1.07}, we can obtain
\begin{eqnarray}
&&\lim_{t\rightarrow0}e^{it\partial_{x}^{2}}f=f,\label{1.08}
\end{eqnarray}
where $f\in H^{s}(\R)$, $s\geq\frac{1}{4}$.
Recently, Bez et al. \cite{BLN2020} extended \eqref{1.07} and obtained a maximal-in-time estimate of the following form
\begin{eqnarray}
&\left\|\sum\limits_{j=1}^{+\infty}\lambda_{j}\left|e^{it\partial_{x}^{2}}f_{j}\right|^{2}
\right\|_{L_{x}^{2,\infty}L_{t}^{\infty}(\mathbf{R}\times I)}\leq C
\|\lambda\|_{\ell^{\beta}},\label{1.09}
\end{eqnarray}
where $\beta<2$, $\lambda=(\lambda_{j})_{j}\in \ell^{\beta}$
and  $(f_j)_{j=1}^{\infty}$ is an orthonormal system in $\dot{H}^{\frac{1}{4}}(\mathbf{R})$.
By using \eqref{1.09}, for $\gamma_{0}\in\mathfrak{S}^{\beta}(\dot{H}^{\frac{1}{4}})$, Bez et al. \cite{BLN2020} proved that
\begin{eqnarray}
&\lim\limits_{t\longrightarrow 0}\rho_{\gamma(t)}(x)=
\rho_{\gamma_{0}}(x),\,\quad a.e.x\in\mathbf{R},\label{1.010}
\end{eqnarray}
where
\begin{eqnarray*}
\rho_{\gamma(t)}=\sum\limits_{j=1}^{+\infty}\lambda_{j}\left|e^{it\partial_{x}^{2}}f_{j}\right|^{2},
\rho_{\gamma_{0}}=\sum\limits_{j=1}^{+\infty}\lambda_{j}\left|f_{j}\right|^{2}
\end{eqnarray*}
denote the density function of
  $\gamma(t), \,\gamma_{0}$, respectively. Here $\gamma(t)=\sum\limits_{j=1}^{+\infty}\lambda_{j}|
e^{it\partial_{x}^{2}}f_{j}\rangle\langle e^{it\partial_{x}^{2}}f_{j}|$
is the solution to operator equation
\begin{eqnarray*}
&&i\frac{d\gamma(t)}{dt}=\left[-\partial_{x}^{2}, \gamma\right], \\
&&\gamma(0)=\gamma_{0}=\sum\limits_{j=1}^{+\infty}\lambda_{j}|f_{j}\rangle\langle f_{j}|.
\end{eqnarray*}
Very recently, by using full randomization introduced by  Hadama and Yamamoto \cite{HY}, for $\gamma_{0}\in\mathfrak{S}^{2}(L^{2})$, Yan et al. established the probabilistic convergence of the density function
  for the Sch\"odinger operator \cite{YDHXY2024} and for the Boussinesq operator \cite{YLYL2024}
   , respectively. The result for the Schr\"{o}dinger operator improves \eqref{1.010} in a probabilistic sense.

\subsection{Motivation and contributions}

By  exploiting  the full randomization on manifolds,
Yan et al. established the probabilistic convergence of the density function
  for the Sch\"odinger operator \cite{YDHXY2024} and the Boussinesq operator \cite{YLYL2024},
respectively.
 A natural question is whether these results extend to the Ostrovsky operator? In this paper,
 we answer this question affirmatively. Specifically, by employing the full randomization of
  compact operators \cite{HY}, for $\gamma_{0}\in\mathfrak{S}^{2}(L^{2})$, we prove that for any  $0<\epsilon\leq1$,
\begin{eqnarray*}
&&\lim\limits_{t \rightarrow 0}(\mathbb{P}\otimes \widetilde{\mathbb{P}})
\left(\left\{(\omega , \widetilde{\omega}) \in(\Omega \times \widetilde{\Omega})
||\rho_{\gamma_{0}^{\omega, \tilde{\omega}}}-\rho_{e^{-t(\partial_{x}^{3}+\partial_{x}^{-1})} \gamma_{0}^{\omega, \widetilde{\omega}}
e^{t(\partial_{x}^{3}+\partial_{x}^{-1})}}|>C\alpha_{1}\right\}\right)=0,
\end{eqnarray*}
which is independent of $x$,
where $(f_j)_{j=1}^{\infty}$ is an orthonormal system in $L^{2}(\R)$
 and  $f_{j}^{\omega}$  is  defined as in
 \eqref{1.013},
\begin{eqnarray*}
&&\alpha_{1}=\left(\left\|\gamma_0\right\|_{\mathfrak{S}^2}+1\right)
e \epsilon^{\frac{1}{2}}\left(\epsilon \ln \frac{1}{\epsilon}\right)^{\frac{3}{2}},
\rho_{\gamma_{0}^{\omega, \tilde{\omega}}}=\sum_{j=1}^{\infty} \lambda_j g_{j}^{(2)}(\widetilde{\omega})
\left|f_{j}^{\omega}\right|^2, \\
&&\rho_{e^{-t(\partial_{x}^{3}+\partial_{x}^{-1})} \gamma_{0}^{\omega, \widetilde{\omega}}
e^{t(\partial_{x}^{3}+\partial_{x}^{-1})}}=\sum_{j=1}^{\infty} \lambda_{j} g_{j}^{(2)}(\widetilde{\omega})\left|e^{-t(\partial_{x}^{3}+\partial_{x}^{-1})} f_{j}^{\omega}\right|^{2},
\end{eqnarray*}
$\rho_{\gamma_{0}^{\omega, \tilde{\omega}}},\rho_{e^{-t(\partial_{x}^{3}+\partial_{x}^{-1})}
 \gamma_{0}^{\omega, \widetilde{\omega}}e^{t(\partial_{x}^{3}+\partial_{x}^{-1})}} $
 are the density function
 of $\gamma_{0}^{\omega, \tilde{\omega}},e^{-t(\partial_{x}^{3}+\partial_{x}^{-1})}
 \gamma_{0}^{\omega, \widetilde{\omega}}e^{t(\partial_{x}^{3}+\partial_{x}^{-1})},$ respectively.

\subsection{Introduction to notation and spaces}
We first fix some notation that will be used throughout the paper.
\begin{eqnarray*}
&&U(t)f=e^{-t(\partial_{x}^{3}+\partial_{x}^{-1})}f=\frac{1}{\sqrt{2\pi}}\int_{\SR}e^{ix\xi}e^{it\phi(\xi)}\mathscr{F}_{x}f(\xi)d\xi,
\end{eqnarray*}
where $\phi(\xi)=\xi^{3}+\xi^{-1}$.

\noindent For
$1\leq\alpha<\infty$, the Schatten norm $\|\gamma\|_{\mathfrak{S}^{\alpha}}$ is defined by
\begin{eqnarray*}
&\|\gamma\|_{\mathfrak{S}^{\alpha}}=:\left(\sum\limits_{j=1}^{\infty}
|\lambda_{j}|^{\alpha}\right)^{\frac{1}{\alpha}},
\end{eqnarray*}
with $(\lambda_{j}^{2})_{j}$ denoting the eigenvalues of  $\gamma^{\ast}\gamma$.
The space of operators with finite $\mathfrak{S}^{\alpha}$
norm is the Schatten class.

Special cases include the Hilbert-Schmidt norm ($\alpha= 2$), which for an integral
 operator corresponds to the $L^{2}$ norm of its kernel, and the operator norm  ($\alpha=\infty$), defined as $\|\gamma\|_{\mathfrak{S}^{\infty}}=\|\gamma\|_{L^{2}\longrightarrow L^{2}}.$

\subsection{Initial value randomization}
To prove Theorem 1.1, we introduce a full randomization procedure for compact operators.
 We begin by randomizing a single function on $\R$.

\noindent{\bf Randomization of single function on $\R$.}
The randomization of  single function on $\R$ originates from
\cite{BOP-2015, BOP2015, ZF2011, ZF2012}.
Firstly, we present the Wiener decomposition of  the frequency space.
 We require that $\psi\in C_{c}^{\infty}(\R)$  is a real-valued,
even, non-negative bump function supported on $[-1,1]$ such that
\begin{eqnarray*}
&&\sum\limits_{k\in \z}\psi(\xi-k)=1,
\end{eqnarray*}
For $\forall k\in \Z$,
  $\psi(D-k)f:\R\rightarrow\mathbb{C}$ is defined as follows:
\begin{eqnarray*}
&&(\psi(D-k)f)(x)=\mathscr{F}_{x}^{-1}\big(\psi(\xi-k)\mathscr{F}_{x}f\big)(x)
=\int_{\SR}e^{ix\xi}\psi(\xi-k)\mathscr{F}_{x}fd\xi,\quad x\in \R.
\end{eqnarray*}
  $(\Omega,\mathcal{A}, \mathbb{P})$, $(\widetilde{\Omega},\widetilde{\mathcal{A}},
     \widetilde{\mathbb{P}})$ are  two different probability spaces.
For each $k\in\mathbf{Z}$,  $\{g_{k}^{(1)}(\omega)\}_{k\in \z}$
  and $\{g_{j}^{(2)}(\widetilde{\omega})\}_{j\in \mathbb{N}^{+}}$  are
   two sequence of independent, zero-mean,
 real-valued random variables with probability distributions $\mu_{k}^{1},$ $\mu_{j}^{2}$,
 respectively.
  Assume that $\mu_{k}^{1}$ and $\mu_{j}^{2}$  satisfy
    the following property:
   for $\forall \gamma_{k}, \gamma_{j}\in \R,\, \forall k\in\Z,\, \forall j\in\mathbb{N}^{+}$,
    $\exists c>0$
\begin{eqnarray}
&&\Big|\int_{-\infty}^{+\infty}e^{\gamma_{k} x}d\mu_{k}^{1}(x)\Big|\leq
 e^{c\gamma_{k}^2},\label{1.011}\\
&&\Big|\int_{-\infty}^{+\infty}e^{\gamma_{j} x}d\mu_{j}^{2}(x)\Big|\leq
 e^{c\gamma_{j}^2}.\label{1.012}
\end{eqnarray}
For $j\in {N}^{+}$ and $f_{j}\in L^2(\R)(j\in\mathbb{N}^{+})$,
we define
\begin{eqnarray}
&&f_{j}^{\omega}=\sum_{k\in\z}g_{k}^{(1)}(\omega)\psi(D-k)f_{j}\label{1.013}.
\end{eqnarray}
 $f_{j}^{\omega}$ is called as the randomization of $f_{j}$.
Moreover, we define
\begin{eqnarray*}
\|f\|_{L_{\omega,\widetilde{\omega}}^{p}(\Omega\times\widetilde{\Omega})}=
\left(\int_{\Omega\times\widetilde{\Omega}}|f(\omega,\widetilde{\omega})|^{p}
d(\mathbb{P}\times \widetilde{\mathbb{P}})\right)^{\frac{1}{p}}=\left(\int_{\Omega}\int_{\widetilde{\Omega}}
|f(\omega,\widetilde{\omega})|^{p}d\mathbb{P}(\omega)
d\widetilde{\mathbb{P}}(\widetilde{\omega})
\right)^{\frac{1}{p}}.
\end{eqnarray*}

\noindent{\bf Full randomization of compact operator on $\R$.}
The full randomization of compact operator on $\R$ originated from \cite{HY}. Let $\gamma_{0}$ be a compact operator, which
   admits the following singular-value decomposition
\begin{eqnarray}
&&\gamma_{0}=\sum_{j=1}^{\infty}\lambda_{j}|f_{j}\rangle\langle f_{j}|,\label{1.014}
\end{eqnarray}
where $\lambda_{j}\in\mathbb{C}$,  $(f_j)_{j=1}^{\infty}$ is an orthonormal system in $L^{2}(\R)$.

For any $(\omega, \widetilde{\omega})\in \Omega\times \widetilde{\Omega}$,
 the full randomization of the compact operator $\gamma_{0}=\sum\limits_{j=1}^{\infty}\lambda_{j}
 |f_{j}\rangle\langle f_{j}|$ is defined as follows:
\begin{eqnarray}
&&\gamma_{0}^{\omega,\widetilde{\omega}}=\sum\limits_{j=1}^{\infty}\lambda_{j}g_{j}^{(2)}
(\widetilde{\omega})|f_{j}^{\omega}\rangle\langle f_{j}^{\omega}|,\label{1.015}
\end{eqnarray}
where  $(f_j)_{j=1}^{\infty}$ is an orthonormal system in $L^{2}(\R)$ and $f_{j}^{\omega}$ is the randomization of
  $f_{j}$  which is defined as in \eqref{1.013}.

From \eqref{1.015}, by using  $|u\rangle\langle v| f=\langle v,f\rangle u$ and
the definition of density function of integral operator, we have
\begin{eqnarray}
&\hspace{-1cm}\rho_{\gamma_{0}^{\omega, \tilde{\omega}}}=\sum\limits_{j=1}^{\infty} \lambda_j g_{j}^{(2)}(\widetilde{\omega})
\left|f_{j}^{\omega}\right|^2,\, \rho_{e^{-t(\partial_{x}^{3}+\partial_{x}^{-1})} \gamma_{0}^{\omega, \widetilde{\omega}}e^{t(\partial_{x}^{3}+\partial_{x}^{-1})}}=\sum\limits_{j=1}^{\infty} \lambda_{j} g_{j}^{(2)}(\widetilde{\omega})\left|e^{-t(\partial_{x}^{3}+\partial_{x}^{-1})} f_{j}^{\omega}\right|^{2}.\label{1.016}
\end{eqnarray}

\subsection{Some comments about new ingredient established in this paper}

The main new ingredient is Lemma 2.4. Now we explain the reason.

Due to the singularity of the phase function $\xi^{3}+\frac{1}{\xi}$
at the origin,   we  can not completely  follow the methods in \cite{YDHXY2024}.
We will prove Theorem 1.1 by employing Lemma 2.4, which is derived
  from Lemmas 2.6 and 3.2 in  \cite{YZDY2022} and Lemma 6.3 in \cite{YDHXY2024}.

\subsection{The main result}
\begin{Theorem}\label{Theorem1} (Stochastic continuity at zero related to Schatten
 norm on $\R$)
Let $r \in[2, \infty)$, $\gamma_{0} \in \mathfrak{S}^{2}(L^{2})$, and let $(f_j)_{j=1}^{\infty}$ be an orthonormal system in $L^{2}(\R)$. Then, for every $0<\epsilon\leq1$, there exists
  $0<\delta\leq M_{2}^{-4}\epsilon^{6}$,
   such that for $|t|<\delta$, we have
\begin{eqnarray}
&&\left\|F(t, x, \omega, \widetilde{\omega})\right\|_{L_{\omega,\widetilde{\omega}}^{r}
(\Omega\times\widetilde{\Omega})}
\leq Cr^{\frac{3}{2}}\epsilon^{2}(\|\gamma_{0}\|_{\mathfrak{S}^{2}}+1),\label{1.017}
\end{eqnarray}
where $f_{j}^{\omega}$ is defined as in
 \eqref{1.013}, and the constant $M_{2}$  appears in \eqref{2.016}; the constant $C$ is independent of $r$, $\epsilon$, $x$ and $t$.
Moreover, we have
\begin{eqnarray}
&&\lim\limits_{t \rightarrow 0}(\mathbb{P}\otimes \widetilde{\mathbb{P}})
\left(\left\{(\omega , \widetilde{\omega}) \in(\Omega \times \widetilde{\Omega}) ||F(t, x, \omega, \widetilde{\omega})|>C\alpha_{1}\right\}\right)=0,\label{1.018}
\end{eqnarray}
which is independent of $x$,
where
$F(t, x, \omega, \widetilde{\omega})=\sum\limits_{j=1}^{\infty} \lambda_{j} g_{j}^{(2)}(\widetilde{\omega})\left(|f_{j}^{\omega}|^{2}-|e^{-t(\partial_{x}^{3}
+\partial_{x}^{-1})}f_{j}^{\omega}|^{2}\right),\,\alpha_{1}=\left(\left\|\gamma_0\right\|_{\mathfrak{S}^2}+1\right)
e \epsilon^{\frac{1}{2}}\left(\epsilon \ln \frac{1}{\epsilon}\right)^{\frac{3}{2}}$.

\end{Theorem}
\noindent{\bf Remark 1.} Motivated by the idea of Lemmas 2.6, 3.2 of \cite{YZDY2022}
 and Lemma 6.3 of  \cite{YDHXY2024},  we present the proof of Theorem 1.1. Moreover, Theorem 1.1 extends Theorem 1.3 of Yan et al. \cite{YZDY2022} to the setting of orthonormal systems in $L^{2}(\R)$.

\section{Probabilistic estimates of some random series}
\setcounter{equation}{0}

\setcounter{Theorem}{0}

\setcounter{Lemma}{0}

\setcounter{section}{2}

This section is devoted to  giving probabilistic estimates of some random series.

\begin{Lemma}\label{lem2.1}(Khinchin type inequalities)
We assume that $\{g^{(1)}_ {k}(\omega)\}_{k\in \z}$ and $\{g^{(2)}_{j}
(\widetilde{\omega})\}_{j\in\mathbb{N}^{+}}$ are independent, zero mean, real valued
 random variable sequences equipped with probability distributions $\mu^{1}_{k}(k\in\Z)$ and
  $\mu^{2}_{j}(j\in\mathbb{N}^{+})$ satisfying   \eqref {1.011}
    and \eqref {1.012} on probability spaces $(\Omega,\mathcal{A}, \mathbb{P})$ and $(\widetilde{\Omega},\widetilde{\mathcal{A}},\widetilde{\mathbb{P}})$,
    respectively. Then,  for $r\in[2, \infty)$ and $(a_k)_k,\,(b_j)_{j}\in\ell^2$,  we have
\begin{eqnarray}
&&\left\|\sum_{k\in\z}a_{k}g_{k}^{(1)}(\omega)\right\|_{L_{\omega}^{r}(\Omega)}\leq
 C_{k}r^{\frac{1}{2}}\|a_{k}\|_{\ell_{k}^{2}},\label{2.01}\\
&&\left\|\sum_{j=1}^{\infty}b_{j}g_{j}^{(2)}
(\widetilde{\omega})\right\|_{L_{\widetilde{\omega}}^{r}(\widetilde{\Omega})}
\leq C_{j}r^{\frac{1}{2}}\|b_{j}\|_{\ell_{j}^{2}}.\label{2.02}
\end{eqnarray}
Here  $C_k ,\, C_j$ are independent of $r$.

\end{Lemma}

 From \cite[ Lemma 3.1]{BT2008I}, we know that Lemma 2.1 is valid.

\begin{Lemma}\label{lem2.2} Assume that $F(t, x, \omega, \widetilde{\omega}):[-1,1]\times\R
\times \Omega \times \widetilde{\Omega}\longrightarrow \R$
 is  measurable and $r \geq 2$, and that for every $0<\epsilon\leq1$, there exists $\delta>0$, such that for all $|t|<\delta$,
 we have
\begin{eqnarray}
&\|F(t, x,\omega, \widetilde{\omega})\|_{L_{\omega, \tilde{\omega}}^{r}(\Omega \times
\tilde{\Omega})} \leq C r^{\frac{3}{2}} \epsilon^{2}\left(\left\|\gamma_0\right\|_{\mathfrak{S}^2}+1\right),\label{2.03}
\end{eqnarray}
where the constant $C$ is independent of $r$, $\epsilon$, $x$ and $t$.
Then, we have
\begin{eqnarray}
\lim _{t \rightarrow 0}(\mathbb{P} \otimes \widetilde{\mathbb{P}})
\left(\left\{(\omega , \widetilde{\omega}) \in\Omega \times \widetilde{\Omega}||F(t, x,\omega, \widetilde{\omega})|>C\alpha_{1}
\right\}\right)=0,\label{2.04}
\end{eqnarray}
which is independent of $x$, where $\alpha_{1}=\left(\left\|\gamma_{0}\right\|_{\mathfrak{S}^2}+1\right) e \epsilon^{\frac{1}{2}}\left(\epsilon \ln \frac{1}{\epsilon}\right)^{\frac{3}{2}}$.
\end{Lemma}

From \cite[Lemma 6.6]{YDHXY2024},  we know that Lemma 2.2 is valid.

\noindent{\bf Remark 2.} In this paper, we apply Lemma 2.2 with
\begin{eqnarray*}
&&F(t, x, \omega, \widetilde{\omega})=\sum\limits_{j=1}^{\infty} \lambda_{j} g_{j}^{(2)}(\widetilde{\omega})\left(|f_{j}^{\omega}|^{2}-|e^{-t(\partial_{x}^{3}
+\partial_{x}^{-1})}f_{j}^{\omega}|^{2}\right),
\end{eqnarray*}
 and Lemma 2.2 plays an important role in the proof of \eqref{1.018} in Theorem 1.1.

\begin{Lemma}\label{lem2.3}
Assume that $|k|\leq M(M\geq1)$, and $\|f_{j}\|_{L^{2}}=1$. Then, for all $0<\epsilon\leq1$, we have
\begin{eqnarray}
&&\sum\limits_{|k|\leq M}
\left|\int_{\SR} \psi(\xi-k)\left(e^{it\phi(\xi)}-1\right)e^{i x\xi}
\mathscr{F}_x f_{j}(\xi) d \xi\right|^{2}\nonumber\\
&&\leq C\left(\epsilon^{4}
+M^{3}|t|^{2}\epsilon^{-8}+M^{8}|t|^{2}\right).\label{2.05}
\end{eqnarray}
Here $\phi(\xi)=\xi^{3}+\xi^{-1}$, the constant $C$ is independent of $M$, $\epsilon$, $x$, $t$. In particular, there exists $\delta>0$ with $\delta\leq M^{-4}\epsilon^{6}$ such that whenever
 $|t|<\delta$, we have
\begin{eqnarray}
&&\left(\sum_{|k| \leq M}\left|\int_{\SR} \psi(\xi-k)
\left(e^{it\phi(\xi)}-1\right) e^{i x\xi}\mathscr{F}_x f_j(\xi)
d \xi\right|^2\right)^{\frac{1}{2}} \leq C\epsilon^{2}.\label{2.06}
\end{eqnarray}

\end{Lemma}

\begin{proof} Inspired by \cite{YZDY2022}, we prove Lemma 2.3.
By using the Cauchy-Schwartz inequality and  $|\psi(\xi-k)\left(e^{it\phi(\xi)}-1\right)|\leq C$, we have
\begin{eqnarray}
&&\int_{|\xi|\leq \frac{\epsilon^{4}}{M}} \left|\psi(\xi-k)\left(e^{it\phi(\xi)}-1\right)
e^{i x\xi}\mathscr{F}_x f_{j}(\xi)\right| d \xi\nonumber\\
&&\leq C\int_{|\xi|\leq \frac{\epsilon^{4}}{M}}|\mathscr{F}_{x}f_{j}|d\xi\leq C\left(\int_{|\xi|\leq \frac{\epsilon^{4}}{M}}
d\xi\right)^{\frac{1}{2}}\|f_{j}\|_{L^{2}}\leq CM^{-\frac{1}{2}}\epsilon^{2}.\label{2.07}
\end{eqnarray}
For $\frac{\epsilon^{4}}{M}\leq|\xi|\leq 1$, we have
\begin{eqnarray}
&&|e^{it\phi(\xi)}-1|\leq C|t||\xi|^{-1}\leq C|t|M\epsilon^{-4}.\label{2.08}
\end{eqnarray}
For $1\leq|\xi|\leq M+1$, we have
\begin{eqnarray}
&&|e^{it\phi(\xi)}-1|\leq C|t||\xi|^{3}\leq C|t|M^{3}.\label{2.09}
\end{eqnarray}
Combining \eqref{2.07}-\eqref{2.09} with H\"{o}lder inequality  yields
\begin{eqnarray}
&&\sum\limits_{|k|\leq M}
\left|\int_{\SR} \psi(\xi-k)\left(e^{it\phi(\xi)}-1\right)e^{i x\xi}\mathscr{F}_x
f_{j}(\xi) d \xi\right|^{2}\nonumber\\
&&\leq C\sum\limits_{|k|\leq M}
\left(\int_{\SR} \left|\psi(\xi-k)\left(e^{it\phi(\xi)}-1\right)\mathscr{F}_x
f_{j}(\xi)\right| d \xi\right)^{2}\nonumber\\
&&\leq C\sum\limits_{|k|\leq M} \left(\int_{|\xi|\leq \frac{\epsilon^{4}}{M}} \left|\psi(\xi-k)\left(e^{it\phi(\xi)}-1\right)
\mathscr{F}_x f_{j}(\xi)\right| d \xi\right)^{2}\nonumber\\
&&+C\sum\limits_{|k|\leq M}\left(\int_{\frac{\epsilon^{4}}{M}\leq|\xi|\leq 1} \left|\psi(\xi-k)\left(e^{it\phi(\xi)}-1\right)
\mathscr{F}_x f_{j}(\xi)\right| d \xi\right)^{2}\nonumber\\
&&+ C\sum\limits_{|k|\leq M}\left(\int_{1\leq|\xi|\leq M+1} \left|\psi(\xi-k)\left(e^{it\phi(\xi)}-1\right)
\mathscr{F}_x f_{j}(\xi)\right| d \xi\right)^{2}\nonumber\\
&&\leq CMM^{-1}\epsilon^{4}
+C\sum\limits_{|k|\leq M}\left(\int_{\frac{\epsilon^{4}}{M}\leq|\xi|\leq 1} |t|M\epsilon^{-4}
\left|\mathscr{F}_x f_{j}(\xi)\right| d \xi\right)^{2}\nonumber\\
&&+ C\sum\limits_{|k|\leq M}\left(\int_{1\leq|\xi|\leq M+1} |t|M^{3}\left|\mathscr{F}_x f_{j}(\xi)\right| d \xi\right)^{2}\nonumber\\
&&\leq CMM^{-1}\epsilon^{4}
+CMM^{2}|t|^{2}\epsilon^{-8}\|f_{j}\|_{L^{2}}^{2}+CMM^{7}|t|^{2}\|f_{j}\|_{L^{2}}^{2}\nonumber\\
&&\leq C\left(\epsilon^{4}
+M^{3}|t|^{2}\epsilon^{-8}+M^{8}|t|^{2}\right).\label{2.010}
\end{eqnarray}

The proof of Lemma 2.3 is finished.
\end{proof}

\noindent{\bf Remark 3.}
Lemma 2.3 will be used in the low-frequency estimates of Lemma 2.4; see specifically \eqref{2.018} in the proof of Lemma 2.4.

\begin{Lemma}\label{lem2.4}(The estimate related to $\ell_{j}^{2}L_{\omega}^{p}$ on $\R$)
Let  $r\geq2$, $(\lambda_{j})_{j}\in \ell_{j}^{2}$, and $(f_{j})_{j=1}^{\infty}$ be an orthonormal system with $\|f_{j}\|_{L^{2}}=1$. Then, for every $0<\epsilon\leq1$, there exists $\delta>0$ with $\delta\leq M_{2}^{-4}\epsilon^{6}$ such that whenever
 $|t|<\delta$, we derive
\begin{align}
&\|\lambda_{j}\|U(t)f_{j}^{\omega}-f_{j}^{\omega}\|_{L_{\omega}^{r}}\|_{\ell_{j}^{2}}\leq
Cr^{\frac{1}{2}}\epsilon^{2}\left(\|\gamma_{0}\|_{\mathfrak{S}^{2}}+1\right),\label{2.011}
\end{align}
where $\|\gamma_{0}\|_{\mathfrak{S}^{2}}=
\left(\sum\limits_{j=1}^{\infty}|\lambda_{j}|^{2}\right)^{\frac{1}{2}}$, the constant $M_{2}$ will appear in \eqref{2.016}, and $f_{j}^{\omega}$ is
 defined as in \eqref{1.013}; the constant $C$ is independent of $r$, $\epsilon$, $x$ and $t$.
\end{Lemma}

\begin{proof}
By using Lemma 2.1, we obtain
\begin{eqnarray}
&\|\lambda_{j}\|U(t)f_{j}^{\omega}-f_{j}^{\omega}\|_{L_{\omega}^{r}}\|_{\ell_{j}^{2}}\leq Cr^{\frac{1}{2}}\left(\sum\limits_{j=1}^{\infty}|\lambda_{j}|^{2}Q_{j}\right)^{\frac{1}{2}},
\label{2.012}
\end{eqnarray}
where
\begin{eqnarray*}
&&Q_{j}=\sum\limits_{k \in \z}
\left|\int_{\SR} \psi(\xi-k)\left(e^{it\phi(\xi)}-1\right)e^{i x\xi}\mathscr{F}_x f_j(\xi)
d \xi\right|^{2}.
\end{eqnarray*}
From $(\lambda_{j})_{j}\in\ell^{2}$, we know that for $\forall\epsilon>0$,
$\exists M_{1}\in \mathbb{N}^{+}$, the following inequality
\begin{eqnarray}
&\left(\sum\limits_{j=M_{1}+1}^{\infty}|\lambda_{j}|^{2}\right)^{\frac{1}{2}}<
\frac{\epsilon^{2}}{2}\label{2.013}
\end{eqnarray}
is valid.

\noindent By using \eqref{2.013} and the Cauchy-Schwartz inequality, we obtain
\begin{align}
\left(\sum_{j=M_{1}+1}^{\infty}|\lambda_{j}|^{2}Q_{j}\right)^{\frac{1}{2}}
&\leq C\left(\sum_{j=M_{1}+1}^{\infty}|\lambda_{j}|^{2}\sum_{k\in\z}\int_{\SR}|\psi(\xi-k)
\mathscr{F}_{x}f_{j}(\xi)|^{2}d\xi\right)^{\frac{1}{2}}\nonumber\\
&\leq C\left(\sum_{j=M_{1}+1}^{\infty}|\lambda_{j}|^{2}\|f_{j}\|_{L^{2}}^{2}\right)^{\frac{1}{2}}\leq C\left(\sum_{j=M_{1}+1}^{\infty}|\lambda_{j}|^{2}\right)^{\frac{1}{2}}\leq C\epsilon^{2},\label{2.014}
\end{align}
where we used the following fact
$$|e^{it\phi(\xi)}-1|\leq2,\,\,\|f_{j}\|_{L^{2}}=1.$$
We claim that for $1\leq j\leq M_{1}$, and for every $0<\epsilon\leq1$, there exists $\delta>0$ with $\delta\leq M_{2}^{-4}\epsilon^{6}$ such that whenever
 $|t|<\delta$,
\begin{eqnarray}
&Q_{j}^{\frac{1}{2}}\leq C\epsilon^{2}\label{2.015}
\end{eqnarray}
is valid, where $M_2$ depends on $M_{1}$. The idea of \cite[Lemma 6.3]{YDHXY2024} and
 \cite[Lemma 3.2]{YZDY2022} will be  used
 to establish \eqref{2.015}.

 By using  line 14 of 1908 in \cite{YZDY2022}, we obtain
\begin{eqnarray*}
&&\sum_{k\in\z}\|\psi(\xi-k)\mathscr{F}_{x}f_{j}\|_{L^{2}}^{2}\leq \|f_{j}\|_{L^{2}}^{2}\leq 3\sum_{k\in\z}\|\psi(\xi-k)\mathscr{F}_{x}f_{j}\|_{L^{2}}^{2},
\end{eqnarray*}
consequently, for $\forall\epsilon>0$,  $\exists M_{2}\geq 2026$, we have
\begin{eqnarray}
&\left(\sum\limits_{|k|\geq M_{2}}
\|\psi(\xi-k)\mathscr{F}_{x}f_{j}\|_{L^{2}}^{2}\right)^{\frac{1}{2}}\leq
\frac{\epsilon^{2}}{2}.\label{2.016}
\end{eqnarray}
Since $\supp \psi\subset[-1,1]$ and $|e^{it\phi(\xi)}-1|\leq2$, by using \eqref{2.016}
 and the Cauchy-Schwartz inequality, we obtain
\begin{eqnarray}
&&\left(\sum\limits_{|k| \geq M_{2}}\left|\int_{\SR} \psi(\xi-k)\left(e^{it\phi(\xi)}-1\right)
 e^{i x\xi} \mathscr{F}_x f_j(\xi) d \xi\right|^2\right)^{\frac{1}{2}} \nonumber\\
&&\leq 2\left(\sum\limits_{|k| \geq M_2} \int_{\mathbf{R}}\left|\psi(\xi-k)
\mathscr{F}_x f_j(\xi)\right|^2 d \xi\right)^{\frac{1}{2}}\nonumber\\
&&=2\left(\sum\limits_{|k| \geq M_2}\left\|\psi(\xi-k)
\mathscr{F}_x f_j\right\|_{L^2}^2\right)^{\frac{1}{2}}\leq \epsilon^{2}.\label{2.017}
\end{eqnarray}
By using Lemma 2.3, for every $0<\epsilon\leq1$, there exists $\delta>0$ with $\delta\leq M_{2}^{-4}\epsilon^{6}$ such that whenever
 $|t|<\delta$, we have
\begin{eqnarray}
&&\left(\sum_{|k| \leq M_2}\left|\int_{\SR} \psi(\xi-k)
\left(e^{it\phi(\xi)}-1\right) e^{i x\xi}\mathscr{F}_x f_j(\xi)
d \xi\right|^2\right)^{\frac{1}{2}} \leq C\epsilon^{2}.\label{2.018}
\end{eqnarray}
From \eqref{2.017} and \eqref{2.018}, we have that
 \eqref{2.015} is valid.

By using \eqref{2.012}, \eqref{2.014} and \eqref{2.015},
we have
\begin{eqnarray}
&&\left\|\lambda_{j}\|U(t)f_{j}^{\omega}-f_{j}^{\omega}\right\|_{L_{\omega}^{p}}\|_{\ell_{j}^{2}}
\leq C r^{\frac{1}{2}}\left(\sum_{j=1}^{+\infty}\left|\lambda_j\right|^2 Q_{j}\right)^{\frac{1}{2}}\nonumber\\
&&\leq C r^{\frac{1}{2}}\left(\sum_{j=1}^{M_{1}}\left|\lambda_j\right|^2Q_{j}\right)^{\frac{1}{2}}
+C r^{\frac{1}{2}}\left(\sum_{j=M_{1}+1}^{+\infty}\left|\lambda_j\right|^2 Q_{j}\right)^{\frac{1}{2}}\nonumber\\
&&\leq Cr^{\frac{1}{2}}\sup_{1\leq j\leq M_{1}}Q_{j}^{\frac{1}{2}}
\left(\sum_{j=1}^{M_{1}} |\lambda_j|^2\right)^{\frac{1}{2}}+C r^{\frac{1}{2}}\epsilon^{2}\nonumber\\
&&\leq C r^{\frac{1}{2}}\left(\epsilon^{2}\left(\sum_{j=1}^{\infty} |\lambda_j|^2\right)^{\frac{1}{2}}+\epsilon^{2}\right)
\leq C r^{\frac{1}{2}} \epsilon^{2}\left(\left\|\gamma_0\right\|_{\mathfrak{S}^2}+1\right).
\label{2.019}
\end{eqnarray}

The proof of Lemma 2.4 is finished.

\end{proof}

\noindent{\bf Remark 4.}
Lemma 2.4 will be used in the proof of \eqref{1.017} in Theorem 1.1; specifically, see the last inequality in \eqref{3.02} within the proof of Lemma 3.1.

\begin{Lemma}\label{lem2.5}
Let $r\geq2$, and let $(f_j)_{j=1}^{\infty}$ be an orthonormal system in $L^{2}(\R)$. Define
 $f_{j}^{\omega}$ as in \eqref{1.013}.
 Then, we have
\begin{eqnarray}
&&\left\|f_{j}^{\omega}\right\|_{L_{\omega}^{r}} \leq C r^{\frac{1}{2}}\left\|f_{j}\right\|_{L^{2}}
=C r^{\frac{1}{2}},\quad\left\|U(t) f_{j}^{\omega}\right\|_{L_{\omega}^{r}} \leq C r^{\frac{1}{2}}\left\|f_{j}\right\|_{L^{2}}
=C r^{\frac{1}{2}},\label{2.020}
\end{eqnarray}
where the constant $C$ is independent of $r$, $t$, $x$.
\end{Lemma}

By using a proof  similar to \cite[Lemma 2.3]{BOP2015}, we obtain  Lemma 2.5.

\noindent{\bf Remark 5.}
Lemma 2.5 will be used in the proof of \eqref{1.017} in Theorem 1.1; specifically, see the second-to-last inequality in \eqref{3.02} in the proof of Lemma 3.1.

\bigskip
\bigskip

\section{Proof of Theorem  1.1: Stochastic continuity at zero related
 to Schatten norm on $\R$}

\setcounter{equation}{0}

\setcounter{Theorem}{0}

\setcounter{Lemma}{0}

\setcounter{section}{3}

In this section, we present the stochastic continuity at zero related to
 Schatten norm on $\R$.

Inspired by \cite{YDHXY2024}, we give the following estimates that  play
 a key role in proving Theorem 1.1.

\begin{Lemma}\label{lem3.1}
Let $r\geq2$, $(\lambda_{j})_{j}\in \ell_{j}^{2}$,
and let $(f_j)_{j=1}^{\infty}$ be an orthonormal system in $L^{2}(\R)$.
Then, for every $0<\epsilon\leq 1$, there exists $\delta>0$ with $\delta\leq M_{2}^{-4}\epsilon^{6}$ such that whenever $|t|<\delta$, we have
\begin{align}
&&\left\|\sum\limits_{j=1}^{\infty} \lambda_{j} g_{j}^{(2)}(\widetilde{\omega})
 I_{j}^{(1)}\right\|_{L_{\omega, \tilde{\omega}}^{r}(\Omega \times \tilde{\Omega})}
\leq Cr^{\frac{3}{2}}\epsilon^{2}\left(\|\gamma_{0}\|_{\mathfrak{S}^{2}}+1\right),\label{3.01}
\end{align}
where
\begin{eqnarray*}
&&I_{j}^{(1)}=\left(e^{-t(\partial_{x}^{3}+\partial_{x}^{-1})} f_{j}^{\omega}-f_{j}^{\omega}\right) e^{-t(\partial_{x}^{3}+\partial_{x}^{-1})}\overline{f}_{j}^{\omega},\,  \|\gamma_{0}\|_{\mathfrak{S}^{2}}=\left(\sum\limits_{j=1}^{\infty}|\lambda_{j}|^{2}\right)^{\frac{1}{2}},
\end{eqnarray*}
the constant $M_{2}$ appears in \eqref{2.016},  and $f_{j}^{\omega}$ denotes the randomization of $f_{j}$ defined in \eqref{1.013}; the constant $C$ is independent of $r$, $\epsilon$, $x$ and $t$.
\end{Lemma}

\begin{proof}
Combining Minkowski inequality, H\"older inequality with Lemmas 2.1, 2.4 and 2.5, we find that
for any $0<\epsilon\leq 1$, there exists $\delta>0$ with
 $0<\delta\leq M_{2}^{-4}\epsilon^{6}$ such that
   whenever $|t|<\delta$, we have
\begin{align}
\left\|\sum_{j=1}^{\infty} \lambda_{j} g_{j}^{(2)}(\widetilde{\omega})
I_{j}^{(1)}\right\|_{L_{\omega, \tilde{\omega}}^{r}(\Omega \times \tilde{\Omega})}
&\leq C r^{\frac{1}{2}}\left\|\lambda_{j}  I_{j}^{(1)}\right\|_{L_\omega^{r} \ell_{j}^{2}}
\leq C r^{\frac{1}{2}}\left\|\lambda_{j}  I_{j}^{(1)}\right\|_{\ell_{j}^{2} L_{\omega}^{r}}\nonumber\\
&\leq C r^{\frac{1}{2}}\left\|\lambda_{j}\left\|\left(e^{-t(\partial_{x}^{3}+\partial_{x}^{-1})} f_{j}^{\omega}-f_{j}^{\omega}\right)\right\|_{L_{\omega}^{2r}}\left\| e^{-t(\partial_{x}^{3}+\partial_{x}^{-1})}\overline{f}_{j}^{\omega}\right\|_{L_{\omega}^{2r}}\right\|_{\ell_{j}^{2}}\nonumber\\
&\leq C r\left\|\lambda_{j}\left\|\left(e^{-t(\partial_{x}^{3}+\partial_{x}^{-1})} f_{j}^{\omega}-f_{j}^{\omega}\right)\right\|_{L_{\omega}^{2r}}\right\|_{\ell_{j}^{2}}\nonumber\\
&\leq C r^{\frac{3}{2}} \epsilon^{2}\left(\left\|\gamma_{0}\right\|_{\mathfrak{S}^2}+1\right).\label{3.02}
\end{align}

The proof of Lemma 3.1 is finished.
\end{proof}

\begin{Lemma}\label{lem3.2}
Let $r\geq2$, $(\lambda_{j})_{j}\in \ell_{j}^{2}$, and let $(f_j)_{j=1}^{\infty}$ be an orthonormal system in $L^{2}(\R)$.
Then, for every $0\epsilon\leq1$, there exists $\delta>0$ with $\delta\leq M_{2}^{-4}\epsilon^{6}$ such that whenever $|t|<\delta$, we have
\begin{align}
&&\left\|\sum\limits_{n=1}^{\infty} \lambda_{j} g_{j}^{(2)}(\widetilde{\omega})
I_{j}^{(2)}\right\|_{L_{\omega, \tilde{\omega}}^{r}(\Omega \times \tilde{\Omega})}
\leq Cr^{\frac{3}{2}}\epsilon^{2}\left(\|\gamma_{0}\|_{\mathfrak{S}^{2}}+1\right),\label{3.03}
\end{align}
where
\begin{eqnarray*}
&&I_{j}^{(2)}=f_{j}^{\omega}\left(e^{-t(\partial_{x}^{3}+\partial_{x}^{-1})}
\overline{f}_{j}^{\omega}-\overline{f}_{j}^{\omega}\right),\,
\|\gamma_{0}\|_{\mathfrak{S}^{2}}=\left(\sum\limits_{j=1}^{\infty}
|\lambda_{j}|^{2}\right)^{\frac{1}{2}},
\end{eqnarray*}
the constant $M_{2}$
 appears in \eqref{2.016}, and $f_{j}^{\omega}$ denotes the randomization of $f_{j}$ defined in \eqref{1.013}; the constant $C$ is independent of $r$, $\epsilon$, $x$ and $t$.
\end{Lemma}

\begin{proof}
By applying  the method similar to that used in the proof of Lemma 3.1,
we can see that Lemma 3.2 holds.

The proof of Lemma 3.2 is finished.
\end{proof}

\noindent{\bf Proof of Theorem  1.1}. Now, we use Lemmas 3.1,
 3.2 and 2.2 to prove Theorem 1.1. Obviously, we have
\begin{eqnarray}
&&\left|e^{-t(\partial_{x}^{3}+\partial_{x}^{-1})} f_{j}^{\omega}\right|^2
-\left|f_{j}^{\omega}\right|^2=I_{j}^{(1)}+I_{j}^{(2)}.\label{3.05}
\end{eqnarray}
By using \eqref{3.05} and  Lemmas 3.1 and 3.2, we have
\begin{eqnarray}
&&\left\|\sum_{j=1}^{\infty} \lambda_{j} g_{j}^{(2)}(\widetilde{\omega})
\left|f_{j}^{\omega}\right|^{2}-\sum_{j=1}^{\infty} \lambda_{j}
g_{j}^{(2)}(\widetilde{\omega})\left|e^{-t(\partial_{x}^{3}
+\partial_{x}^{-1})}f_{j}^{\omega}\right|^2\right\|_{L_{\omega,
\tilde{\omega}}^{r}(\Omega \times \widetilde{\Omega})}\nonumber\\
&&=\left\|\sum_{j=1}^{\infty} \lambda_{j} g_{j}^{(2)}(\widetilde{\omega})
\left(I_{j}^{(1)}+I_{j}^{(2)}\right)\right\|_{L_{\omega, \tilde{\omega}}^{r}
(\Omega \times \widetilde{\Omega})}\nonumber\\
&&\leq\left\|\sum_{j=1}^{\infty} \lambda_{j} g_{j}^{(2)}(\widetilde{\omega})
 I_{j}^{(1)}\right\|_{L_{\omega, \tilde{\omega}}^{r}(\Omega \times \tilde{\Omega})}
 +\left\|\sum_{j=1}^{\infty} \lambda_{j} g_{j}^{(2)}(\widetilde{\omega})
 I_{j}^{(2)}\right\|_{L_{\omega, \tilde{\omega}}^{r}(\Omega \times \tilde{\Omega})}\nonumber\\
&&\leq C r^{\frac{3}{2}} \epsilon^{2}
\left(\left\|\gamma_{0}\right\|_{\mathfrak{S}^2}+1\right).\label{3.06}
\end{eqnarray}
From\eqref{3.06}, we see that \eqref{1.017} is valid. Combining \eqref{1.017}
 and Lemma 2.2, we obtain that \eqref{1.018} holds.

The proof of Theorem 1.1 is finished.

\bigskip











\end{document}